\documentclass[a4paper,11pt]{article}
\usepackage{a4wide}
\usepackage{theorem}
\usepackage{amsmath}
\usepackage{array}
\usepackage{amssymb}
\usepackage{amsfonts}
\usepackage[french,english]{babel}
\usepackage{epsf}
\usepackage{epsfig}

\newtheorem{theo}{\indent Theorem\newline}[section]

{\theorembodyfont{\rmfamily}%
\newtheorem{rem}[theo]{\noindent Remark}}
{\theorembodyfont{\rmfamily}%
 \theoremstyle{break}%
}
\newtheorem{prop}[theo]{\indent Proposition\newline}
\newtheorem{lemma}[theo]{\indent Lemma\newline}
\newtheorem{cor}[theo]{\indent Corollary\newline}

 \def\N{{\mathbb{N}}}
\def\Z{{\mathbb{Z}}}

\def\R{{\mathbb{R}}}
\def\C{{\mathbb{C}}}

\usepackage{eufrak}
\newcommand{\goth}[1]{\EuFrak{#1}}

\newcommand{\p}{\mathop{\goth p}\nolimits}

\newcommand{\conj}{\mathop{\rm conj}\nolimits}
\newcommand{\coker}{\mathop{\rm coker}\nolimits}

\newlength{\indentation}%
\makeatletter
\newcommand\@makefntextsans[1]{%
    \parindent 0em%
    \noindent%
    \hb@xt@0em{\hss}%
    #1}
\def\footnotetextsans{%
     \@ifnextchar [\@xfootnotenextsans%
       {\@footnotetextsans}}
\def\@xfootnotenextsans[#1]{%
  \begingroup%
     \csname c@\@mpfn\endcsname #1\relax%
  \endgroup%
  \@footnotetextsans}
\long\def\@footnotetextsans#1{\insert\footins{%
    \reset@font\footnotesize%
    \interlinepenalty\interfootnotelinepenalty%
    \splittopskip\footnotesep%
    \splitmaxdepth \dp\strutbox \floatingpenalty \@MM%
    \hsize\columnwidth \@parboxrestore%
    \color@begingroup%
      \@makefntextsans{%
        \rule\z@\footnotesep\ignorespaces#1\@finalstrut\strutbox}
    \color@endgroup}}
\makeatother

\begin{document}

\cleardoublepage
\title{Enumerative invariants of stongly semipositive real symplectic six-manifolds}
\author{Jean-Yves Welschinger}
\maketitle

\makeatletter\renewcommand{\@makefnmark}{}\makeatother
\footnotetextsans{Keywords : Real symplectic manifold, 
rational curve, enumerative geometry.}
\footnotetextsans{AMS Classification : 53D45, 14N35, 14N10, 14P99.}

{\bf Abstract:}

Following the approach of Gromov and Witten \cite{Gro,Wit}, we define invariants under deformation
of stongly semipositive real symplectic six-manifolds. These invariants provide lower bounds in real enumerative geometry, namely for
the number of real rational $J$-holomorphic curves which realize a given homology class
and pass through a given real configuration of points. 

\section*{Introduction}

A smooth compact symplectic manifold $(X, \omega)$ of dimension $2n$ is said to be {\it semipositive} 
as soon as for every spherical class $d \in H_2 (X ; \Z)$ such that $[\omega] d > 0$, the implication
$c_1(X) d \geq 3 - n \implies c_1(X) d \geq 0$ holds. These manifolds provide a favourable framework
to define genus zero Gromov-Witten invariants (see \cite{McDSal}) and in particular these
invariants are enumerative. We will assume throughout the paper that the manifold $(X, \omega)$ is 
{\it strongly semipositive}, by which we mean that for every spherical class $d \in H_2 (X ; \Z)$ 
such that $[\omega] d > 0$, the implication $c_1(X) d \geq 2 - n \implies c_1(X) d \geq 1$ holds
(see Remark \ref{remcomment}).
An important source of examples are smooth projective Fano manifolds. The manifold is said to be
{\it real} when it is equipped with an antisymplectic involution, that is an involution $c_X$ such
that $c_X^* \omega = - \omega$. We denote by $\R X$ the fixed locus of this involution and
by $\R {\cal J}_\omega$ the set of almost complex structures of $X$ of finite regularity $C^l$, 
$l \gg 1$, which are tamed by $\omega$ and for which $c_X$ is $J$-antiholomorphic. In this context,
one can count the number of real rational $J$-holomorphic curves which realize a given homology class
$d \in H_2 (X ; \Z)$ and pass through some given real configuration of points, under the assumptions
that $J$ is generic and this number finite. However, this number does depend in general on the choice 
of the almost complex structure $J$ or the configuration of points, 
basically
because the field of real numbers is not algebraically closed. In \cite{Wels0}, \cite{Wels1}, a way 
of counting 
these real rational $J$-holomorphic curves with respect to some sign $\pm 1$ has been introduced in 
order to define an integer which neither depends on the choice of $J \in \R {\cal J}_\omega$ nor on the
choice of the points. This integer is invariant under deformation of the real symplectic manifold, but
has been defined only in dimension four. The integer $\pm 1$ depended on the parity of the number of 
real isolated double points of the real rational $J$-holomorphic curves which were counted. These
curves are indeed singular in general in this dimension.
The question appeared then whether it was possible to obtain similar results in higher dimensional 
real symplectic manifolds. A partial answer to this question was 
obtained in \cite{Wels2} where a way to count such real rational curves in real algebraic convex
$3$-manifolds has been introduced. The sign $\pm 1$ depended then on some spinor state of the real 
rational curve which
was first defined and which required the choice of a $Pin^-_3$ structure on the real locus $\R X$.
The integers thus defined were invariants under isomorphism of real algebraic convex
$3$-manifolds. However, as was pointed out in \cite{Wels2}, very few such convex
$3$-manifolds indeed have non trivial genus zero Gromov-Witten invariants, namely $\C P^3$, 
$\C P^2 \times \C P^1$, $(\C P^1)^3$,
$Fl(\C^3)$ and the quadric in $\C P^4$ for the ones I am aware of. The aim of this work is to
build such integer valued invariants for any strongly semipositive real symplectic six-manifolds, 
see Theorem \ref{theosemipos}. These integers are obtained by counting the number of real 
rational $J$-holomorphic curves which realize a given homology class
and pass through a given real configuration of points with respect to some spinor state in
$\{ \pm 1 \}$ that we first define, see \S \S \ref{subsubsectspinorstate} and \ref{subsectGCR}. 
They are invariant under strongly semipositive
deformation of the real symplectic six-manifold. This means that if $\omega_t$ is 
a continuous family of strongly semipositive symplectic forms on $X$ for which $c_X^* \omega_t = 
-\omega_t$, then these invariants are the same for all triples $(X , \omega_t , c_X)$. Moreover,
they provide lower bounds in real enumerative geometry, namely for
the number of real rational $J$-holomorphic curves which realize a given homology class
and pass through a given real configuration of points, see Corollary \ref{corlowerbounds}. 
In a first part of this paper, we
define spinor states for real rational curves in  real algebraic convex manifolds of any dimension.
In order to define the spinor states of such real rational $J$-holomorphic curves, it is necessary 
here to make some topological assumption on the real locus
of the manifolds, basically that its second Stiefel-Whitney class either vanishes or equals the square
of the first Stiefel-Whitney class, see \S \ref{subsubsectHYP} for the exact hypothesis and Remark
\ref{remcomment} for a comment on these hypothesis. 
We then extend this definition to real rational curves of strongly semipositive real symplectic manifolds
and define the invariants in dimension six.
%

\tableofcontents

\section{Preliminaries}

\subsection{Moduli space of genus zero stable maps}
\label{subsectstable}

Let $(X , c_X)$ be a smooth real algebraic convex manifold of complex dimension $n \geq 3$, that 
is a smooth real 
projective manifold such that for every morphism $u : \C P^1 \to X$, the vanishing $H^1 (\C P^1 ;
u^* TX) = 0$ occurs. Denote by $\R X = \text{fix} (c_X)$ its real locus, which we assume to be 
nonempty. Let $d \in H_2 (X ; \Z)$ be such that $(c_X)_* d = -d$ and $(n-1) / (c_1(X)d - 2)$. We set
$k_d = \frac{1}{n-1} (c_1(X)d - 2) + 1$. Note that $k_d \in \N^*$ as soon as $d$ is realized by some
rational curve, see Lemma $11$ of \cite{FP}. Denote by $\overline{\cal M}^d_{k_d} (X)$ the space of
genus zero stable maps of $X$ which realize $d$ and have $k_d$ marked points. Let $ev^d :
\overline{\cal M}^d_{k_d} (X) \to X^{k_d}$ be the evaluation map. Note that 
$\dim_\C \overline{\cal M}^d_{k_d} (X) = c_1(X)d + n - 3 + k_d = nk_d$ so that $ev^d$ is a morphism
between projective manifolds of the same dimension. Let $\tau \in \sigma_{k_d}$ be such that
$\tau^2 = id$. Following $\S 1.1$ of \cite{Wels2}, we denote by $c_\tau : (x_1 , \dots , x_{k_d})
\in X^{k_d} \mapsto \big( c_X ( x_{\tau (1)}) , \dots , c_X ( x_{\tau (k_d)}) \big)$ the
associated real structure on $X^{k_d}$. In the same way, denote by $c_{\overline{\cal M} , \tau}$
the real structure of $\overline{\cal M}^d_{k_d} (X)$ induced by $c_{{\cal M} , \tau} :
(u , z_1 , \dots , z_{k_d}) \in {\cal M}or_d (X) \times (\C P^1)^{k_d} \mapsto (c_X \circ u \circ
\conj , \conj(z_{\tau (1)}) , \dots , \conj(z_{\tau (k_d)})) \in {\cal M}or_d (X) \times 
(\C P^1)^{k_d}$, where $\conj$ is the standard complex conjugation of $\C P^1$ and
${\cal M}or_d (X) = \{ u : \C P^1 \to X \, | \, u_* [\C P^1] = d \}$ 
(see Theorem $1.1$ of \cite{Wels2}). Denote by $\R_\tau X^{k_d} = \text{fix} (c_\tau)$ and
$\R_\tau \overline{\cal M}^d_{k_d} (X) = \text{fix} (c_{\overline{\cal M} , \tau})$ the real loci
of these spaces. The evaluation morphism restricts to $\R_\tau ev^d : \R_\tau 
\overline{\cal M}^d_{k_d} (X) \to \R_\tau X^{k_d}$. Finally, denote by 
$\overline{\cal M}^d_{k_d} (X)^*$ (resp. $\R_\tau \overline{\cal M}^d_{k_d} (X)^*$) the subspace of
simple stable maps of $\overline{\cal M}^d_{k_d} (X)$ (resp. $\R_\tau \overline{\cal M}^d_{k_d} (X)$),
it is contained in the smooth locus of $\overline{\cal M}^d_{k_d} (X)$, see \cite{FP}. Note that
the singular locus of $\overline{\cal M}^d_{k_d} (X)$ is of codimension greater that one and will
play no r\^ole in this paper.
\begin{theo}
\label{theoprelim}
Let $(X , c_X)$ be a smooth real algebraic convex manifold of complex dimension $n \geq 3$.

1) The divisor $\text{Red} = \overline{\cal M}^d_{k_d} (X)^* \setminus {\cal M}^d_{k_d} (X)^*$ has 
normal crossings.

2) Let $(u , C , \underline{z}) \in {\cal M}^d_{k_d} (X)^*$ and ${\cal N}_u$ be its normal sheaf.
Then the isomorphisms $\ker d|_{(u , C , \underline{z})} ev^d \cong H^0 (C ;  {\cal N}_u \otimes
{\cal O}_C (-\underline{z})) = H^0 (C ;  N_u \otimes
{\cal O}_C (-\underline{z})) \oplus H^0 (C ;  {\cal N}_u^{\text{sing}} \otimes
{\cal O}_C (-\underline{z}))$ and $\coker d|_{(u , C , \underline{z})} ev^d \cong 
H^1 (C ;  {\cal N}_u \otimes {\cal O}_C (-\underline{z}))$ hold.

3) As soon as $(n , c_1 (X) d) \neq (3,4)$, the locus of stable maps $(u , C , \underline{z}) \in 
{\cal M}^d_{k_d} (X)^*$ for which $u$ is not an immersion is mapped onto some submanifold of $X^{k_d}$
having codimension greater than one. 
\end{theo}
\begin{rem}
In the third part of this theorem, the condition $n \geq 3$ is crucial. For a discussion on the
condition $(n , c_1 (X) d) \neq (3,4)$, see Remark $3.2$ of \cite{Wels2}.
\end{rem}
{\bf Proof:}

The first part is Theorem $3$ of \cite{FP}. The proof of the second part goes exactly along the same
lines as the one of Lemma $1.3$ of \cite{Wels2}, it is not reproduced here. Now the locus of
stable maps $(u , C , \underline{z}) \in {\cal M}^d_{k_d} (X)^*$ for which $u$ is not an immersion
and for which $u^* TX \otimes {\cal O}_C (-1)$ is generated by its global sections is of codimension 
$n-1$ in ${\cal M}^d_{k_d} (X)^*$. This follows from the fact that the tautological
section $\sigma : (u , C , z) \in {\cal M}^d_1 (X)^* \mapsto d|_z u \in \text{Hom}_\C 
(T_z C , T_{u(z)} X)$ vanishes transversely at those points, which can be proved as Proposition $3.1$ 
of \cite{Wels2}. From the second part of Theorem \ref{theoprelim}, the locus of stable maps 
$(u , C , \underline{z}) \in 
{\cal M}^d_{k_d} (X)^*$ for which $\dim H^1 (C ;  {\cal N}_u \otimes {\cal O}_C (-\underline{z})) > 1$
is mapped onto some submanifold of $X^{k_d}$ having codimension greater than one.
Let then $(u , C , \underline{z}) \in 
{\cal M}^d_{k_d} (X)^*$ be such that $\dim H^1 (C ;  {\cal N}_u \otimes {\cal O}_C (-\underline{z})) 
\leq 1$. From a theorem of Grothendieck (\cite{Grot}), the normal bundle $N_u$ of $u$ is 
isomorphic to ${\cal O}_C (a_1) \oplus \dots \oplus {\cal O}_C (a_{n-1})$, where $a_i \geq k_d - 1$
for $1 \leq i \leq n-2$ and $a_{n-1} \geq  k_d - 2$. If $u$ is not an immersion, then
$\deg (N_u) \leq c_1 (X) d - 3 = (n-1)(k_d - 1) - 1$, so that $N_u$ has to be isomorphic to
${\cal O}_C (k_d - 1)^{n-2} \oplus {\cal O}_C (k_d - 2)$. Remember that the convexity of $X$ forces
$u^* TX$ to be a direct sum of line bundles of non-negative degrees, see Lemma $10$ of \cite{FP}.
Assume that $u^* TX \otimes {\cal O}_C (-1)$ is not generated by its global sections, then
$u^* TX \cong F \oplus {\cal O}_C^k$, where $F$ is a direct sum of line bundles of positive degrees
and $k \geq 1$. Let $z_C \in C$ be the point where $du$ vanishes. The latter maps $TC \otimes
{\cal O}_C (z_C)$ to $F$ so that $N_u$ is isomorphic to ${\cal O}_C^k \oplus \big( F/du(TC \otimes
{\cal O}_C (z_C)) \big)$. From what precedes, this forces $k_d \leq 2$ and the inequality 
$\deg (u^* TX)
\geq 3$ forces $k_d \geq 2$. Hence, $k_d = 2$ and $k = 1$. Now the composition of the tautological
section $\sigma : (u , C , z) \in 
{\cal M}^d_1 (X)^* \mapsto d|_z u \in \text{Hom}_\C (T_z C , T_{u(z)} X)$ with the projection
$\text{Hom}_\C (T_z C , T_{u(z)} X) \to \text{Hom}_\C (T_z C , F)$ vanishes transversely since
$F \otimes {\cal O}_C (-1)$ is generated by its global sections. Since $\dim F = n-2$, the result
follows in all the cases but $n=3$, $k_d = 2$. $\square$
\begin{lemma}
\label{lemmarealpart}
Assume that either $n$ is even, or $\tau$ has a fixed point in $\{ 1 , \dots , k_d \}$. Then, as
soon as $\underline{x} \in \R_\tau X^{k_d}$ is generic enough, the fibre $(\R_\tau ev^d)^{-1} 
(\underline{x})$ uniquely consists of irreducible real rational curves having non-empty real parts.
\end{lemma}
Note that when $k_d$ is odd, $\tau$ has a fixed point so that the hypothesis of Lemma 
\ref{lemmarealpart} is satisfied. We will assume throughout the paper that this hypothesis holds.\\

{\bf Proof:}

To begin with, assume that $\tau$ has a fixed point, say $1 \in \{ 1 , \dots , k_d \}$ and denote
$\underline{x}$ by $(x_1 , \dots , x_{k_d})$. Let $(u , C , \underline{z}) \in (\R_\tau ev^d)^{-1} 
(\underline{x})$. From Theorem \ref{theoprelim}, as
soon as $\underline{x}$ is generic enough, $C$ is irreducible. From the definition of 
$c_{{\cal M} , \tau}$, there exists $\phi \in \text{Aut} (C)$ such that $c_X \circ u \circ \conj
= u \circ \phi^{-1}$ and $\conj (z_{\tau (i)}) = \phi (z_i)$, $i \in \{ 1 , \dots , k_d \}$.
Let $c_C = \phi^{-1} \circ \conj$, then $c_X \circ u = u \circ c_C$ and $c_C (z_{\tau (i)}) = z_i$,
$i \in \{ 1 , \dots , k_d \}$. In particular, $c_C$ is a real structure on $C$ which has one fixed 
point at least, namely $x_1$. Assume now that $n$ is even. From what has been done, we can also
assume that $k_d$ is even. From the definition of $k_d$, it implies that $c_1 (X) d$ is odd.
Now if $C$ has an empty real part, then $c_1 (X) d = w_2 (u^* TX) [C] \mod (2) = 
2 w_2 (u^* TX/c_X) [C/c_C] = 0 \mod (2)$, hence the contradiction. $\square$
\begin{lemma}
\label{lemmarm*}
Let $\R {\cal M}^*$ be a connected component of $\R_\tau {\cal M}^d_{k_d} (X)^*$. Then, the homology
class $u_* [\R C] \in H_1 (\R X ; \Z/2\Z)$ does not depend on the choice of $(u, C , \underline{z})
\in \R {\cal M}^*$.
\end{lemma}
The homology class given by Lemma \ref{lemmarm*} will be denoted by $d_{\R {\cal M}^*} 
\in H_1 (\R X ; \Z/2\Z)$.

{\bf Proof:}

Let $\R_\tau U^d \to \R_\tau {\cal M}^*$ be the universal curve. Then, all the fibres
of $\R_\tau U^d$ have same homology class in $H_1 (\R_\tau U^d ; \Z/2\Z)$. The result is thus
obtained after composition with the morphism $H_1 (\R_\tau U^d ; \Z/2\Z) \to H_1 (\R X ; \Z/2\Z)$
induced by the evaluation map $\R_\tau U^d \to \R X$. $\square$

\subsection{Spinor states}
\label{subsectspinor}

\subsubsection{$\widetilde{GL}_m^{\pm} (\R)$-structures}
\label{subsubsectHYP}

Denote by $\widetilde{GL}_m (\R)$ the universal covering of $GL_m (\R)$. It can be equipped with
two different group structures which turn the covering map into a morphism. Denote by 
$\widetilde{GL}_m^+ (\R)$ (resp. $\widetilde{GL}_m^- (\R)$) the group structure for which the lift
of a reflexion is of order two (resp. of order four). Let $M \to \R X$ be a vector bundle of rank $m$
and $R_M$ be the associated $GL_m (\R)$-principal bundle of frames. The obstruction to the
existence of a $\widetilde{GL}_m^+ (\R)$-principal bundle (resp. $\widetilde{GL}_m^- (\R)$-principal 
bundle) $P_M$ over $R_M$ is carried by the characteristic class $w_2 (M) \in H^2 (\R X ;
\Z/2\Z)$. (resp. $w_2 (M) + w_1^2 (M) \in H^2 (\R X ; \Z/2\Z)$), see \cite{KT} for example. From now 
on, we will assume that one of the following holds and will denote by {\bf HYP} these hyposthesis.

1) If either $k_d$ is odd or $n=3 \mod (4)$, we assume that $0 \in \{ w_2 (\R X) , w_2 (\R X)  
+ w_1^2 (\R X) \}$. We then set $M = T \R X$ and equip this bundle with a $\widetilde{GL}_m^+ (\R)$
or $\widetilde{GL}_m^- (\R)$-structure depending on whether $w_2 (\R X)$ or $w_2 (\R X)  
+ w_1^2 (\R X)$ vanishes.

2) If $k_d$ is even and $n=0 \mod (4)$, we assume that $w_2 (\R X) = 0$.  
We then set $M = T \R X \oplus \text{Det} (\R X)^3$, where $\text{Det} (\R X)$ is the determinant line
bundle of $\R X$, and equip this bundle with a spin structure,
which is possible since $w_1 (M) = w_2 (M) = 0$.

3) If $k_d$ is even and $n=2 \mod (4)$, we assume that $w_2 (\R X) =w_1^2 (\R X)$. 
We then set $M = T \R X \oplus \text{Det} (\R X)$ and we equip this bundle with a spin structure,
which is possible since $w_1 (M) = w_2 (M) = 0$.

4) If $k_d$ is even and $n=1 \mod (4)$, we assume that there exists $w \in H^1 (\R X ; \Z/2\Z)$
such that $w^2 \in \{ w_2 (\R X) , w_2 (\R X)  + w_1^2 (\R X) \}$. We then set $M = T \R X \oplus 
L_{\R X} (w)^2$, where $L_{\R X} (w)$ is the line bundle over $\R X$ whose Euler class is $w$. 
Moreover, we equip this bundle with a $\widetilde{GL}_m^+ (\R)$
or $\widetilde{GL}_m^- (\R)$-structure depending on whether $w^2 = w_2 (\R X)$ or $w^2 = w_2 (\R X)  
+ w_1^2 (\R X)$. Note that $w_1 (M) = w_1 (\R X)$ and $w_2 (M) = w_2 (\R X)  + w^2$.

\begin{rem}
\label{remw}
Under the assumption $4$, we will have actually to restrict ourselves to real rational curves $A$
whose real locus $\R A$ satisfy $<w , [\R A]> \neq 0$. Using the terminology of Lemma \ref{lemmarm*},
this means that we will restrict ourselves to connected components $\R {\cal M}^*$ of 
$\R_\tau {\cal M}^d_{k_d} (X)^*$ for which $<w , d_{\R {\cal M}^*}> \neq 0 \in \Z/2\Z$.
\end{rem}

{\bf Examples:}

1) If $n=3$, then from Wu formula, $w_2 (\R X) =w_1^2 (\R X)$ so that condition $1$ is always
satisfied, see page $132$ of \cite{Miln} for instance.

2) All projective spaces satisfy these hypothesis. Indeed, if $X = \C P^n$ and $w$ is the generator
of $H^1 (\R P^n ; \Z/2\Z)$, then $w_1 (\R X) = (n+1) w$ and $w_2 (\R X) = \frac{n(n+1)}{2} w^2$. If
$n \neq 1 \mod (4)$ or $k_d$ is odd, one easily check the hypothesis. If $k_d$ is even and 
$n=1 \mod (4)$, one has to choose $w$ to be the generator of $H^1 (\R P^n ; \Z/2\Z)$. It turns out
that the definition of $k_d$ forces the degrees of rational curves to be odd in this case, so that
their real parts are always non trivial against $w$.

3) If $\R X$ is spin and $n \neq 1 \mod (4)$, HYP are satisfied.

\subsubsection{Spinor states of balanced real rational curves}
\label{subsubsectspinorstate}

Let $u : (\C P^1 , \conj) \to (X , c_X)$ be a $\Z / 2\Z$-equivariant immersion and $0 \to T \C P^1
\to u^* TX \to N_u \to 0$ be the associated exact sequence. Let $d = u_* [\C P^1]$ and $k_d$ be the
associated integer, see \S \ref{subsectstable}. We assume that the hypothesis HYP of \S 
\ref{subsubsectHYP} are satisfied and deduce from the previous exact sequence a splitting
$u^* T \R X = T \R P^1 \oplus \R N_u$ which is well defined up to homotopy. Assume now that
the vanishing $H^1 (\C P^1 ; N_u \otimes {\cal O}_{\C P^1} (-k_d)) = 0$ occurs, so that
$N_u$ gets isomorphic to ${\cal O}_{\C P^1} (k_d - 1)^{n-1}$. Such a curve is said to be 
{\it balanced}, as its normal bundle is a direct sum of isomorphic lines bundles. Let $\xi \in \R P^1$
and $(v_1 , \dots , v_{n-1})$ be a basis of the fibre $\R N_u|_\xi$. For $i \in \{ 1 , \dots ,
n-1 \}$, there exists a unique holomorphic line sub bundle of degree $k_d - 1$ of $N_u$ which
contains $v_i$. This immediatly follows from the isomorphism $P (N_u) \cong \C P^1 \times \C P^{n-2}$
and the fact that such line sub bundles are mapped onto constant sections of 
$\C P^1 \times \C P^{n-2}$. We deduce a decomposition $\R N_u = L_{\R P^1} (k_d - 1)^{n-1}$ which
is well defined up to homotopy.

Now let $M$ be the vector bundle defined in \S \ref{subsubsectHYP}. In all cases $1$ to $4$
considered in \S \ref{subsubsectHYP}, we obtained a decomposition $u^* M = L_{\R P^1} (0) \oplus
L_{\R P^1} (k_d - 1)^{m-1}$ well defined up to homotopy. That is $u^* M$ is given a splitting as 
the direct sum of $m$ orientable real line bundles when $k_d$ is odd and as the direct sum
of $m-1$ non orientable real line bundles and one orientable one when $k_d$ is even. Note that in
all the cases $u^* M$ is orientable, since $m = 3 \mod (4)$. When both $k_d$ and $n$ are even,
$u^* M$ is even oriented since by assumption $M$ is. When $k_d$ is even and $n$ is odd, we equip
$u^* M$ with an orientation. We are then ready to define the spinor state of the balanced
real rational curve $Im (u)$.

{\bf $1$st case : $k_d$ is odd.} In this case, we choose some trivialization of each factor
$L_{\R P^1} (0)$ of $u^* M = L_{\R P^1} (0)^m$ given by some non vanishing section $v_i$,
$i \in \{ 1 , \dots m \}$. We hence get a loop $(v_1 (\xi), \dots , v_m (\xi))$, $\xi \in \R P^1$,
of the $GL_m (\R)$-bundle $u^* R_M$ of frames of $u^* M$. We define $sp_{\p} (u) = +1$ (resp.
$sp_{\p} (u) = -1$) if this loop does lift (resp. does not lift) to a loop of the
$\widetilde{GL}_m^\pm (\R)$-principal bundle $u^* P_M$ given in \S \ref{subsubsectHYP} which defines 
the $\widetilde{GL}_m^\pm (\R)$-structure $\p$ on $M$. This integer neither depends on the choice of 
the decomposition $u^* M = L_{\R P^1} (0)^m$ nor on the one of the sections 
$v_1 (\xi), \dots , v_m (\xi)$. It also does not depend on the parameterization $u$ of $Im (u)$ and 
is called the {\it spinor state} of $Im (u)$.

{\bf $2$nd case : $k_d$ is even.} Then $m-1 = 2 \mod(4)$ and $u^* M$ is oriented. We choose some 
trivialization of each factor of the decomposition
$u^* M =  L_{\R P^1} (0) \oplus L_{\R P^1} (k_d - 1)^{m-1}$ over $[0 , \pi ] \subset \R P^1$ given 
by some non vanishing section $v_i|_{[0 , \pi ]}$, $i \in \{ 1 , \dots m \}$. We can assume that
$v_1 \in L_{\R P^1} (0)$ and the basis $(v_1 (\xi), \dots , v_m (\xi))$ to be direct. Also, we
assume that the orientations of $\R P^1$ given by $v_1$ and $[0 , \pi ]$ are the same. Then we
can apply the path
$\begin{array}{rcccl}
g(\xi) & = & \left[
\begin{array}{cccc}
1 &&& \\
& \text{rot} (\xi) &&\\
&& \ddots & \\
&&& \text{rot} (\xi) \\
\end{array}
\right] & \in & GL_m (\R)
\end{array}$
to this basis, where $\xi \in [0 , \pi ]$ and 
$\begin{array}{rcl}
\text{rot} (\xi) &=&
\left[
\begin{array}{cc}
\cos (\xi) & -\sin (\xi) \\
\sin (\xi) & \cos (\xi) \\
\end{array}
\right]. \\
\end{array}$ The number of such rotation blocks is odd. We hence get a new decomposition $u^* M = 
L_{\R P^1} (0)^m$ as a direct sum of orientable real line bundles.
\begin{lemma}
\label{lemmahomot}
The homotopy class of the decomposition $u^* M = L_{\R P^1} (0)^m$ hence obtained does not depend 
on the choice of the base $(v_1, \dots , v_m)$.
\end{lemma}
Note that it however depends on the choice of the orientation on $u^* M$.\\

{\bf Proof:}

When $v_1$ is fixed, the orientation of $u^* M$ induces an orientation of the factor 
$L_{\R P^1} (k_d - 1)^{m-1}$. The result follows from the fact that the set of direct basis
of $L_{\R P^1} (k_d - 1)^{m-1}$ is connected. Consider now the base $(-v_1, v_3 , v_2 , \dots , 
v_{m-1} , v_m)$. Since $m = 3 \mod(4)$, it is a direct basis of $u^* M$. Applying the path
$g(\xi)$ taking into account that the orientation of $\R P^1$ is reversed amounts the same as
applying the path $g(\pi - \xi)$ with the same orientation. To conclude, one just has to check that
applying the rotation block $\text{rot} (\xi)$ to the basis $(v_i , v_{i+1})$ has the same effect
as applying the rotation block $\text{rot} (\pi - \xi)$ to the basis $(v_{i+1} , v_i)$. $\square$\\

Thanks to Lemma \ref{lemmahomot}, we can now define the {\it spinor state} 
$sp_{\p , {\goth o}} (u) \in
\{ \pm 1 \}$ as in the first case, where ${\goth o}$ is the chosen orientation of $u^* M$. This 
integer 
neither depends on the choice of the decomposition
$u^* M = L_{\R P^1} (0)^m$ nor on the one of the sections $v_1 (\xi), \dots , v_m (\xi)$. It also
does not depend on the parameterization $u$ of $Im (u)$. Note that the definition of spinor state
given here extends the one of \S $2.2$ of \cite{Wels2}. Remember that the set of 
$\widetilde{GL}_m^{\pm} (\R)$-structures of $M$ is a principal space over $H^1 (\R X ; \Z / 2\Z)$, 
see page $184$ of \cite{KT} for example. We denote this action by $(w , \p) \mapsto  w.\p$, where
$w \in H^1 (\R X ; \Z / 2\Z)$.
\begin{lemma}
\label{lemmaorient}
Let $u : (\C P^1 , \conj) \to (X , c_X)$ be a $\Z / 2\Z$-equivariant balanced immersion  and
$w \in H^1 (\R X ; \Z / 2\Z)$. Then, $sp_{\p , -{\goth o}} (u) = - sp_{\p , {\goth o}} (u)$ and
$sp_{w.\p , {\goth o}} (u) = (-1)^{<w , u_* [\R P^1]>} sp_{\p , {\goth o}} (u)$.
\end{lemma}

{\bf Proof:}

The composition of the path $g(\xi)$, $\xi \in [0 , \pi]$, applied to the basis $(v_1, \dots , v_m)$
and the path $g(\pi - \xi)$, $\xi \in [0 , \pi]$, applied to the basis $(v_1, v_3 , v_2 , v_4 , \dots 
, v_{m-1} , v_m)$ has a non-vanishing homotopy class in $\pi_1 (GL_m (\R)) \cong \Z/2\Z$.
The relation $sp_{\p , -{\goth o}} (u) = - sp_{\p , {\goth o}} (u)$ follows. The relation 
$sp_{w.\p , {\goth o}} (u) = 
(-1)^{<w , u_* [\R P^1]>} sp_{\p , {\goth o}} (u)$ immediatly follows from the definition of the action
of $H^1 (\R X ; \Z / 2\Z)$ on the set of $\widetilde{GL}_m^{\pm} (\R)$-structures of $M$. $\square$

\section{Moduli space  of real rational pseudo-holomorphic curves}

We recall here the construction of the universal moduli space $\R {\cal M}^d (\underline{x})$ of 
real rational pseudo-holomorphic curves which realize a given homology class $d$ and pass through a 
given configuration of points $\underline{x}$, see \cite{McDSal} and \cite{Wels1}
for the real case in dimension $4$.

Let $d \in H_2 (X ; \Z)$ be such that $(c_X)_* d = -d$ and $(n-1) / (c_1(X)d - 2)$. Let 
$k_d = \frac{1}{n-1} (c_1(X)d - 2) + 1$ and $\underline{x} = (x_1 , \dots , 
x_{k_d}) \in X^{k_d}$ be a real configuration of $k_d$ distinct points.
Denote by $\tau \in \sigma_{k_d}$ the permutation of $\{ 1 , \dots , k_d \}$ induced by $c_X$ 
and assume that $\tau (1) = 1$ if $n$ is odd, see Lemma \ref{lemmarealpart}. Let $S$ be
an oriented $2$-sphere and $\underline{z} = (z_1 , \dots , 
z_{k_d})$ be a set of $k_d$ distinct marked points on it. Let ${\cal J}_S$ be the space of complex 
structures of class $C^l$ on $S$
which are compatible with its orientation, where $l$ is a large enough integer. Likewise, let
${\cal J}_\omega$ be the space of almost complex structures of class $C^l$ of $X$, tamed by $\omega$.
Denote by
$${\cal P} (\underline{x}) = \{ (u, J_S , J) \in L^{k,p} (S , X) \times {\cal J}_S \times
{\cal J}_\omega \, | \, u_* [S] = d \, , \, u (\underline{z}) = \underline{x} \text{ and } 
\sigma_{\overline{\partial}} (u) = 0 \},$$
where $\sigma_{\overline{\partial}} (u) = du + J \circ du \circ J_S$ is the Cauchy-Riemann section
of the $C^{l-k}$ Banach bundle ${\cal E}$ over $L^{k,p} (S , X) \times {\cal J}_S \times
{\cal J}_\omega$ whose fibre over $(u, J_S , J)$ is the separable Banach space $L^{k - 1,p} (S , 
\Lambda^{0,1} S \otimes E_u)$, where $E_u = u^* TX$. Denote by ${\cal P}^* (\underline{x})
\subset {\cal P} (\underline{x})$ the subspace of non-multiple pseudo-holomorphic maps.
We recall the following Proposition (see Proposition $3.2.1$ of \cite{McDSal}).
\begin{prop}
The space ${\cal P}^* (x)$ is a separable Banach manifold of class $C^{l-k}$ whose 
tangent
space at $(u,J_S , J) \in {\cal P}^* (x)$ is the space $T_{(u,J_S , J)} {\cal P}^* 
(x) =
\{ (v , \dot{J}_S , \dot{J}) \in L^{k,p} (S , E_u) \times T_{J_S} {\cal J}_S \times T_J 
{\cal J}_\omega
\, | \, v(\underline{z}) = 0 \text{ and } \nabla_{(v , \dot{J}_S , \dot{J})} 
\sigma_{\overline{\partial}} 
 = 0 \}$. $\square$
\end{prop}
Let us fix some $c_X$-invariant riemannian metric $g$ on $X$ and denote by $\nabla$ the connection
induced on $TX$ as well as on all the assiciated vector bundles. Then, 
$\nabla_{(v , \dot{J}_S , \dot{J})} \sigma_{\overline{\partial}}$ writes
$\nabla \sigma_{\overline{\partial}} (v , \dot{J}_S , \dot{J}) =
Dv + J \circ du \circ \dot{J}_S + \dot{J} \circ du \circ J_S$,
where $D$ is the Gromov operator defined by
$v \in  L^{k,p} (S , E_u)  \mapsto D(v) = \nabla v + J \circ \nabla v  
\circ J_S +
\nabla_v J \circ du \circ J_S \in {\cal E}|_{(u,J_S , J)}$.

Denote by ${\cal D}iff (S,z)$ the group of diffeomorphisms of class $C^{l+1}$ of $S$,
which either preserve the orientation and fix $\underline{z}$, or reverse the orientation
and induce the permutation on $\underline{z}$ associated to $\tau$.
Let ${\cal D}iff^+ (S,z)$ (resp. ${\cal D}iff^- (S,z)$) be the subgroup 
of ${\cal D}iff (S,z)$ of orientation
preserving diffeomorphisms (resp. its complement in ${\cal D}iff (S,z)$) and
$s_*$ be the morphism ${\cal D}iff (S,z) \to
\Z / 2\Z$ of kernel ${\cal D}iff^+ (S,z)$.
The group ${\cal D}iff (S,z)$ acts on ${\cal P}^* (x)$ by 
$$\phi . (u,J_S , J) = \left\{ \begin{array}{rcl}
(u \circ \phi^{-1} , (\phi^{-1})^*J_S , J) & \text{if} & s_*(\phi)=+1, \\
(c_X \circ u \circ \phi^{-1} , (\phi^{-1})^*J_S , \overline{c_X}^*(J)) & 
\text{if} & s_*(\phi)=-1,
\end{array} \right. $$
where $(\phi^{-1})^* J_S = s_*(\phi)
d\phi \circ J_S \circ d\phi^{-1}$ and $\overline{c_X}^*(J) = - d c_X \circ J \circ d c_X$. 
The order two elements of 
${\cal D}iff^- (S,z)$ are the only ones in ${\cal D}iff (S,z) \setminus \{ id \}$ which may have 
a non-empty fixed point set in ${\cal P}^* (x)$. In particular, two such involutions have disjoint 
fixed point sets, compare Lemma $1.3$ of \cite{Wels1}. Let $c_S \in {\cal D}iff (S,z)$ be such an 
element,
we denote by $\R {\cal P}^* (x)_{c_S}$ its fixed locus in ${\cal P}^* (x)$. Denote by $\R {\cal J}_S$
(resp. $\R {\cal J}_\omega$) the fixed locus of $c_S$ (resp. $c_X$) in ${\cal J}_S$
(resp. ${\cal J}_\omega$). Likewise, denote by $L^{k,p} (S , E_u)_{+1} = \{ v \in L^{k,p} (S , E_u)
\, | \, d c_X \circ v \circ c_S = v \}$ the fixed locus of $c_S$ in $L^{k,p} (S , E_u)$. Then,
$T_{(u,J_S , J)} \R {\cal P}^* (x)_{c_S} =
\{ (v , \dot{J}_S , \dot{J}) \in L^{k,p} (S , E_u)_{+1} \times T_{J_S} \R {\cal J}_S \times T_J 
\R {\cal J}_\omega \, | \, v(\underline{z}) = 0 \text{ and } \nabla_{(v , \dot{J}_S , \dot{J})} 
\sigma_{\overline{\partial}}  = 0 \}$, see Proposition $1.4$ of \cite{Wels1}. Note that
$\sigma_{\overline{\partial}}$ and $D$ are ${\cal D}iff (S,z)$-equivariant so that $D$ induces
an operator $D_\R : L^{k,p} (S , E_u)_{+1} \to L^{k-1,p} (S , \Lambda^{0,1} S \otimes E_u)_{+1} =
\{ \alpha \in L^{k-1,p} (S , \Lambda^{0,1} S \otimes E_u) \,  | \, d c_X \circ  \alpha \circ c_S =
\alpha \}$. Denote by $L^{k,p} (S , E_{u , - \underline{z}})_{+1} = \{ v \in L^{k,p} (S , E_u)_{+1}
\,  | \, v(\underline{z}) = 0 \}$ and by $H^0_D (S , E_{u , - \underline{z}})_{+1}$ (resp.
$H^1_D (S , E_{u , - \underline{z}})_{+1}$) the kernel (resp. cokernel) of the operator
$D_\R : L^{k,p} (S , E_{u , - \underline{z}})_{+1} \to  
L^{k-1,p} (S , \Lambda^{0,1} S \otimes E_u)_{+1}$. Remember that this operator $D_\R$ induces
a quotient operator $\overline{D}_\R : L^{k,p} (S , E_{u,-\underline{z}})_{+1} / 
du (L^{k,p} (S , TS_{-\underline{z}})_{+1}) \to L^{k-1,p} (S, \Lambda^{0,1} S 
\otimes_\C E_{u})_{+1} / du (L^{k-1,p} (S, \Lambda^{0,1} S \otimes_\C TS)_{+1})$ (see 
formula $1.5.1$ of \cite{IShev} or \S $1.4$ of \cite{Wels1}). Denote by 
$H^0_D (S , {\cal N}_{u , - \underline{z}})_{+1}$ (resp.
$H^1_D (S , {\cal N}_{u , - \underline{z}})_{+1}$) the kernel (resp. cokernel) of the operator
$\overline{D}_\R$. Finally, the action of ${\cal D}iff^+ (S,z)$ on ${\cal P}^* (x)$ is proper,
fixed point free and with closed complements. 
Denote by ${\cal M}^d (x)$ the quotient of ${\cal P}^* (x)$ by the action of ${\cal D}iff^+ (S,z)$.
The projection $\pi : (u, J_S , J) \in {\cal P}^* (x) \mapsto J \in {\cal J}_\omega$
induces on the quotient a projection ${\cal M}^d (x) \to {\cal J}_\omega$ still denoted
by $\pi$. The manifold ${\cal M}^d (x)$ is equipped with an action of the group 
${\cal D}iff (S,\underline{z}) /{\cal D}iff^+ (S,\underline{z}) \cong
\Z / 2 \Z$ which turns $\pi$ into a $\Z / 2 \Z$-equivariant map. Let us denote by 
$\R {\cal M}^d (x)$ the fixed locus of this action and by $\pi_\R$ 
the induced projection $\R {\cal M}^d (x) \to \R {\cal J}_\omega$.
\begin{prop}
\label{proppir}
The space $\R {\cal M}^d (x)$ is a separable Banach manifold of class $C^{l-k}$ and
$\pi_\R : \R {\cal M}^d (x) \to \R {\cal J}_\omega$ is Fredholm
of vanishing index. Moreover, at $[u, J_S , J] \in
\R {\cal M}^d (x)$, the kernel of $\pi_\R$ is isomorphic to $H^0_D 
(S , {\cal N}_{u,-\underline{z}})_{+ 1}$ and its cokernel to $H^1_D 
(S , N_{u,-\underline{z}})_{+ 1}$. $\square$
\end{prop}
The proof of this proposition is the same as the one of Proposition $1.9$ of \cite{Wels1}. It is
not reproduced here.

\section{Spinor states of generalized real Cauchy-Riemann operators}

Let $(u,J_S , J) \in \R {\cal P}^* (x)_{c_S}$ be such that $u$ is an immersion. The normal bundle
$N_u = u^* TX / du (TS)$ is equipped with a $\Z / 2 \Z$-equivariant operator
$D^N : L^{k,p} (S , N_{u , - \underline{z}}) = \{ v \in L^{k,p} (S , N_u) \, | \, 
v(\underline{z}) = 0 \} \to  L^{k-1,p} (S , \Lambda^{0,1} S \otimes N_u)$. Remember that its
complex linear part is some Cauchy-Riemann operator denoted by $\overline{\partial}$ whereas its
complex antilinear part is some order zero operator denoted by $R (v) = N_J (v , du(.))$ where
$N_J$ is the Nijenhuis tensor of $J$, see  Lemma $1.3.1$ of \cite{IShev}. Such an operator is
called a {\it generalized real Cauchy-Riemann operator}, it is Fredholm of vanisning index.
Denote by $c_N$ the complex antilinear involutive morphism induced by $c_X$ on the complex vector
bundle $N_u$. Denote by $\R N_u$ (resp. $\R S$) the fixed locus of $c_N$ (resp. $c_S$). This vector
bundle $\R N_u$ over $\R S$ has rank $n-1$ and the riemannian metric $g$ on $\R X$ induces a
splitting $u^* T \R X = T \R S \oplus \R N_u$. We assume that the hypothesis HYP of
\S \ref{subsubsectHYP} hold and equip the associated bundle $M \to \R X$ with a $\widetilde{GL}_m^\pm
(\R)$-structure. The aim of this paragraph is to define a spinor state for such 
generalized real Cauchy-Riemann operators when they are isomorphisms. We will begin with
standard real Cauchy-Riemann operators in \S \ref{subsectCR} and then extend to generalized ones in 
\S \ref{subsectGCR}.

\subsection{Spinor states of real Cauchy-Riemann operators}
\label{subsectCR}

Denote by ${\cal O}p_{\overline{\partial}} (N_u)$ the space of Cauchy-Riemann operators of class
$C^{l-1}$ on $N_u$,
it is an affine Banach space spanned by $\Gamma^{l-1} (S , \Lambda^{0,1} S \otimes End_\C (N_u))$,
see Appendix $C1$ of \cite{McDSal}. Denote by $\R {\cal O}p_{\overline{\partial}} (N_u)$ the
sub space of ${\cal O}p_{\overline{\partial}} (N_u)$ made of operators which are $\Z/2\Z$-equivariant
with respect to the actions of $c_N$. For $D \in {\cal O}p_{\overline{\partial}}$, we denote by
$D_{\underline{z}}$ its restriction to $L^{k,p} (S , N_{u , - \underline{z}})$ so that 
$D_{\underline{z}}$ is Fredholm with vanishing index.
\begin{lemma}
The set of operator $D \in \R {\cal O}p_{\overline{\partial}} (N_u)$ for which $D_{\underline{z}}$
is an isomorphism is dense open in $\R {\cal O}p_{\overline{\partial}} (N_u)$. The set of these
operators for which $D_{\underline{z}}$ has a one dimensional cokernel is a one dimensional 
submanifold of $\R {\cal O}p_{\overline{\partial}} (N_u)$. Finally, the set of these
operators for which $D_{\underline{z}}$ has a cokernel of dimension greater than one is included
in a countable union of strata of codimensions greater than one in 
$\R {\cal O}p_{\overline{\partial}} (N_u)$. $\square$
\end{lemma}
Remember that the space ${\cal O}p_{\overline{\partial}} (N_u)$ 
(resp. $\R {\cal O}p_{\overline{\partial}} (N_u)$) corresponds to the space of holomorphic structures
on $N_u$ (resp. for which $c_N$ is antiholomorphic), see \cite{Kob}. Let 
$D \in \R {\cal O}p_{\overline{\partial}} (N_u)$ such that $D_{\underline{z}}$
is an isomorphism. Such an operator is said to be {\it balanced}, since it defines a 
holomorphic structure on $N_u$ such that the 
isomorphism $N_u \cong {\cal O}_S (k_d - 1)^{n-1}$ holds. We can thus define as in \S 
\ref{subsubsectspinorstate} the {\it spinor state} $sp_{\p , {\goth o}} (D) \in
\{ \pm 1 \}$ of the balanced operator $D$.
\begin{prop}
\label{propwall}
Let $D^1 , D^2 \in \R {\cal O}p_{\overline{\partial}} (N_u)$ be two operators belonging to two
adjacent connected components of balanced operators of $\R {\cal O}p_{\overline{\partial}} (N_u)$.
Then, $sp_{\p , {\goth o}} (D^1) = - sp_{\p , {\goth o}} (D^2)$.
\end{prop} 

{\bf Proof:}

Let $D_0 \in  \R {\cal O}p_{\overline{\partial}} (N_u)$ be an operator in the wall between the
two components containing $D_1$ and $D_2$. In particular, $\dim_\C H^1 (S , N_u \otimes {\cal O}_S 
(-\underline{z})) = 1$ when $N_u$ is equipped with the holomorphic
structure induced by $D_0$, so that $N_u \cong {\cal O}_S (k_d - 1)^{n-3} \oplus {\cal O}_S (k_d)
\oplus {\cal O}_S (k_d - 2)$. Let $U_0$, $U_1$ be the standard atlas of $\C P^1 \cong S$. The
holomorphic bundle $N_u$ is obtained by gluing the two charts $U_0 \times \C P^{n-1}$ and
$U_1 \times \C P^{n-1}$ with the help of the gluing map
$$
\begin{array}{rcl}
\Phi_0 : (U_0 \cap U_1)  \times \C P^{n-1} & \to & (U_1 \cap U_0)  \times \C P^{n-1} \\
(\xi , (v_1 , \dots , v_{n-1})) & \mapsto & \Big( \frac{1}{\xi} ,
\left[
\begin{array}{ccc}
(\frac{1}{\xi})^{k_d - 1} \text{Id}_{n-3} & 0 & 0 \\
0 & (\frac{1}{\xi})^{k_d} & 0 \\
0 & 0 & (\frac{1}{\xi})^{k_d - 2} \\
\end{array}
\right]
\left(
\begin{array}{c}
v_1 \\
\vdots\\
v_{n-1}\\
\end{array}
\right)
\end{array}
\Big).
$$
The operator $D_0$ writes in these charts $D_0 (v_1 , \dots , v_{n-1}) =
(\overline{\partial} v_1 , \dots , \overline{\partial} v_{n-1})$. Let $f : U_0 \to \C$ be such that
$f(\xi) = \frac{1}{\xi}$ if $|\xi| \geq 1$ and $f(\xi) = \overline{\xi}$ if $|\xi| \leq 1 -
\epsilon$. We choose $f$ such that $\underline{z}$ is disjoint from the support of
$\overline{\partial} f$ and such that $\int_{U_0} \overline{\partial} f \wedge d\xi \neq 0$. For
$t \in \R$, we set 
$D_t (v_1 , \dots , v_{n-1}) = (\overline{\partial} v_1 , \dots , \overline{\partial} v_{n-2} ,
\overline{\partial} v_{n-1} + t \overline{\partial} f \otimes v_{n-2})$. Note that
$\frac{d}{dt}|_{t=0} D_t : (v_1 , \dots , v_{n-1}) \mapsto (0 , \dots , 0, \overline{\partial} f 
\otimes v_{n-2})$ induces a non-vanishing morphism $H^0 (S ; N_u \otimes {\cal O}_S (-\underline{z}))
\to  H^1 (S , N_u \otimes {\cal O}_S (-\underline{z}))$, since 
$\int_{U_0} \overline{\partial} f \wedge d\xi \neq 0$. Thus, the path $t \in \R \mapsto
D_t \in \R {\cal O}p_{\overline{\partial}} (N_u)$ is transversal to the wall of non-balanced operators
at $t=0$. Let
$$
\begin{array}{rcl}
A_0 : U_0  \times \C P^{n-1} & \to & U_0  \times \C P^{n-1} \\
(\xi , (v_1 , \dots , v_{n-1})) & \mapsto & \Big( \xi ,
\left[
\begin{array}{ccc}
\text{Id}_{n-3} & 0 & 0 \\
0 & 1 & 0 \\
0 & t f(\xi) & 1 \\
\end{array}
\right]
\left(
\begin{array}{c}
v_1 \\
\vdots\\
v_{n-1}\\
\end{array}
\right)
\end{array}
\Big),
$$
and  $A_1 = \text{Id} : (\xi , (v_1 , \dots , v_{n-1})) \in U_1  \times \C P^{n-1} \mapsto 
 (\xi , (v_1 , \dots , v_{n-1})) \in U_1  \times \C P^{n-1}$.
This $0$-cochain provides an isomorphism between $N_u$ equipped with the operator $D_t$ and
the holomorphic vector bundle defined by the transition function
$$
\begin{array}{rcl}
\Phi_t : (U_0 \cap U_1)  \times \C P^{n-1} & \to & (U_1 \cap U_0)  \times \C P^{n-1} \\
(\xi , (v_1 , \dots , v_{n-1})) & \mapsto & \Big( \frac{1}{\xi} ,
\left[
\begin{array}{ccc}
(\frac{1}{\xi})^{k_d - 1} \text{Id}_{n-3} & 0 & 0 \\
0 & (\frac{1}{\xi})^{k_d} & 0 \\
0 & t(\frac{1}{\xi})^{k_d - 1}  & (\frac{1}{\xi})^{k_d - 2} \\
\end{array}
\right]
\left(
\begin{array}{c}
v_1 \\
\vdots\\
v_{n-1}\\
\end{array}
\right)
\end{array}
\Big).
$$
We have to prove that for $t \neq 0$, $sp_{\p , {\goth o}} (D_t) = - sp_{\p , {\goth o}} (D_{-t})$.
For this purpose, we assume that $k_d$ is odd, the case $k_d$ even can treated in the same way. Let
$P : U_0 \to \C$ be a polynomial of degree $k_d - 1$ with real coefficients and no real roots.
From the definition of spinor state, the sections $v_1 = (P , 0 , \dots , 0) , \dots ,
v_{n-2} = (0 , \dots , 0 , \xi P , -tP)$ of $N_u$ restricted to $\R S$ can be used to compute this
spinor state. Indeed, they generate $n-2$ real holomorphic line sub bundles of degree $k_d - 1$
of $N_u$ which are in direct sum and provide a trivialization of the real loci of these bundles. 
The choice of the
$(n-1)^{th}$ such line sub bundle being unique up to homotopy, it is not necessary to introduce it.
Now it turns out that as $t$ crosses $0$, the section $v_{n-2}|_{\R S}$ crosses transversely
the zero section of $\R N_u$. The loops of the principal bundle of frames of $\R N_u$ defined
by $(v_1 , \dots , v_{n-2})$ for $t > 0$ and $t < 0$ are thus obtained one from another by
performing a full twist around the axis generated by $v_1 , \dots , v_{n-3}$. The result follows.
$\square$

\subsection{Spinor states of generalized real Cauchy-Riemann operators}
\label{subsectGCR}

Denote now by ${\cal O}p_{\overline{\partial} + R} (N_u)$ the space of generalized Cauchy-Riemann 
operators of class $C^{l-1}$ on $N_u$,
it is an affine Banach space spanned by $\Gamma^{l-1} (S , \Lambda^{0,1} S \otimes End_\R (N_u))$,
see Appendix $C1$ of \cite{McDSal}. Denote by $\R {\cal O}p_{\overline{\partial} + R} (N_u)$ the
sub space of ${\cal O}p_{\overline{\partial} + R} (N_u)$ made of operators which are 
$\Z/2\Z$-equivariant
with respect to the actions of $c_N$. For $D \in {\cal O}p_{\overline{\partial} + R}$, we denote by
$D^{\underline{z}}_\R$ the operator $L^{k,p} (S , N_{u , - \underline{z}})_{+1} \to
L^{k,p} (S , \Lambda^{0,1} S \otimes N_{u})_{+1}$ so that 
$D^{\underline{z}}_\R$ is Fredholm with vanishing index.
\begin{lemma}
\label{lemmawall}
The set of operator $D \in \R {\cal O}p_{\overline{\partial} + R} (N_u)$ for which 
$D^{\underline{z}}_\R$
is an isomorphism is dense open in $\R {\cal O}p_{\overline{\partial} + R} (N_u)$. The set of these
operators for which $D^{\underline{z}}_\R$ has a one dimensional cokernel is a one dimensional 
submanifold of $\R {\cal O}p_{\overline{\partial} + R} (N_u)$. Finally, the set of these
operators for which $D^{\underline{z}}_\R$ has a cokernel of dimension greater than one is included
in a countable union of strata of codimensions greater than one in 
$\R {\cal O}p_{\overline{\partial} + R} (N_u)$. $\square$
\end{lemma}

Let $D \in \R {\cal O}p_{\overline{\partial} + R} (N_u)$ be a {\it balanced operator}, that is
an operator such that $D^{\underline{z}}_\R$ is an isomorphism. Let
$\overline{\partial} \in \R {\cal O}p_{\overline{\partial}} (N_u)$ be a balanced operator
and $\delta : t \in [ 0 , 1 ] \mapsto D_t \in \R {\cal O}p_{\overline{\partial} + R} (N_u)$ be a
generic path joining $\overline{\partial}$ to $D$. Denote by $n_\delta$ the number of times
this path crosses the wall of non balanced operators given by Lemma \ref{lemmawall}, each crossing
being transversal since $\delta$ is generic. Since 
the determinant line bundle $Det (D^{\underline{z}}_\R) = \Lambda^{max} 
\ker (D^{\underline{z}}_\R) \otimes \Lambda^{max} \coker (D^{\underline{z}}_\R)$ is trivial over
the affine space $\R {\cal O}p_{\overline{\partial} + R} (N_u)$, the parity of $n_\delta$ does
not depend on $\delta$, see Proposition $A.2.4$ of \cite{McDSal}.
We can then define the {\it spinor state} of the balanced operator $D$ to be
$sp_{\p , {\goth o}} (D) = (-1)^{n_\delta} sp_{\p , {\goth o}} (\overline{\partial}) \in
\{ \pm 1 \}$, where $sp_{\p , {\goth o}} (\overline{\partial})$ has been defined in \S \ref{subsectCR}.
It follows from Proposition \ref{propwall} that this spinor state $sp_{\p , {\goth o}} (D)$ does not
depend on the choice of $\overline{\partial}$ we made and hence is well defined.

\section{Statement of the results}

\subsection{Statements}

Let $(X , \omega , c_X)$ be a strongly semipositive real symplectic manifold of dimension six.
Let $d \in H_2 (X ; \Z)$ be such that $(c_X)_* d = -d$, $2 / (c_1(X)d - 2)$ and
$c_1 (X)d > 2$. Let $k_d = \frac{1}{2} (c_1(X)d - 2) + 1$ and $\underline{x} = (x_1 , \dots , 
x_{k_d}) \in X^{k_d}$ be a real configuration of $k_d$ distinct points. We assume that 
$\underline{x}$ has at least one real point, see Lemma \ref{lemmarealpart}. 
We label the connected components of $\R X$ by $(\R X)_1 , \dots , (\R X)_N$ and set
$r_i = \# ( \underline{x} \cap (\R X)_i)$, $i \in \{ 1 , \dots , N \}$ as well as
$r = (r_1 , \dots , r_N)$. Let $J \in \R {\cal J}_\omega$ be generic enough, so that there are only
finitely many connected real rational $J$-holomorphic curves which realize the given homology class 
$d$ and pass through $\underline{x}$. Moreover, these curves are all irreducible, smooth and have 
non-empty real part, see Lemma \ref{lemmarealpart}. Let $h \in H_1 (\R X ; \Z/2\Z)$, we denote by 
${\cal R}_d (\underline{x} , J , h)$ the finite set of those curves whose real part realize $h$.
In case $k_d$ is even, we equip $M|_{x_1}$ with an orientation ${\goth o}$. This orientation induces an 
orientation on $M|_{\R A}$ for every $A \in {\cal R}_d (\underline{x} , J , h)$. Note that for
every such curve $A \in {\cal R}_d (\underline{x} , J , h)$, the normal bundle $N_A$ comes 
equipped with a generalized real Cauchy-Riemann operator $D_A$ which is balanced since $J$ is
generic, see Proposition \ref{proppir}. Thus, spinor states of all the
elements $A \in {\cal R}_d (\underline{x} , J, h)$ are well defined, namely as
$sp_{\p , {\goth o}} (D_A)$, see \S \ref{subsectGCR}. We set
$$\chi_r^{d, h \p , {\goth o}} (\underline{x} , J) = \sum_{A \in {\cal R}_d (\underline{x} , J , h)}
sp_{\p , {\goth o}} (A) \in \Z.$$
\begin{theo}
\label{theosemipos}
Let $(X , \omega , c_X)$ be a strongly semipositive real symplectic manifold of dimension six.
Let $d \in H_2 (X ; \Z)$ be such that $(c_X)_* d = -d$, $2/ (c_1(X)d - 2)$ and
$c_1 (X)d > 2$. Let $k_d = \frac{1}{2} (c_1(X)d - 2) + 1$ and $\underline{x} = (x_1 , \dots , 
x_{k_d}) \in X^{k_d}$ be a real configuration of $k_d$ distinct points with at least one real one. 
Label the connected components of $\R X$ by $(\R X)_1 , \dots , (\R X)_N$ and set
$r_i = \# ( \underline{x} \cap (\R X)_i)$, $i \in \{ 1 , \dots , N \}$ as well as
$r = (r_1 , \dots , r_N)$. Let $h \in H_1 (\R X ; \Z/2\Z)$, assume that the hypothesis HYP of 
\S \ref{subsubsectHYP} hold and denote by $M \to \R X$ the given
rank $m$ vector bundle. In case $k_d$ is even, equip $M|_{x_1}$ with an orientation ${\goth o}$. 
Then, the integer 
$\chi_r^{d, h \p , {\goth o}} (\underline{x} , J)$ neither depends on the generic choice of $J \in
\R {\cal J}_\omega$ nor on the choice of $\underline{x}$.
\end{theo}
See Remark \ref{remcomment} below for a comment on the hypothesis of Theorem \ref{theosemipos}.
It follows from this theorem that the integer $\chi_r^{d, h \p , {\goth o}} (\underline{x} , J)$
can be denoted by $\chi_r^{d, h \p , {\goth o}}$ without ambiguity. We can even get rid of the
orientation ${\goth o}$ from the following.
\begin{cor}
Under the hypothesis of Theorem \ref{theosemipos}, assume that $k_d$ is even. 
If the restriction of $M$ over the connected component of $\underline{x} \cap \R X$ in 
$\R X$ is not
orientable, then $\chi_r^{d, h \p , {\goth o}}= 0$. In particular, the genus zero Gromov-Witten 
invariant $GW_0 (X , d , pt, \dots , pt)$ is even. $\square$
\end{cor}
Note that if all the points of $\underline{x} \cap \R X$ are not in the same connected component
of $\R X$, then $\chi_r^{d, h \p , {\goth o}}$ vanishes since the real locus of rational curves is
connected and thus cannot interpolate points in different connected components.
The generating function for this invariant is then some polynomial $\chi^{d, h \p} (T) =
\sum_{r \in \N^N} \chi_r^{d, h \p} T^r \in \Z [T_1 , \dots ,T_N]$, where 
$T^r = T_1^{r_1} \dots T_N^{r_N}$ and we have set $\chi_r^{d, h \p} = 0$ when it is not well defined.
This polynomial has the same parity as $k_d$ and each of its monomials $\chi_r^{d, h \p} T^r$ 
only depends on one indeterminate as we just saw. Theorem \ref{theosemipos} means that the function
$$\chi^{\p} : (d,h) \in H^2 (X ; \Z) \times H_1 (\R X ; \Z/2\Z) \mapsto \chi^{d, h \p} (T) \in \Z [T]$$
only depends on the real symplectic six-manifold $(X , \omega , c_X)$ and is invariant under
strongly semipositive deformation of this real symplectic six-manifold. This means that if $\omega_t$ is 
a continuous family of strongly semipositive symplectic forms for which $c_X^* \omega_t = -\omega_t$,
then this function is the same for all triples $(X , \omega_t , c_X)$. This invariant immediatly
provides the following lower bounds in real enumerative geometry.
\begin{cor}
\label{corlowerbounds}
Under the hypothesis of Theorem \ref{theosemipos},
$$|\chi_r^{d, h \p}| \leq \# {\cal R}_d (\underline{x} , J , h) \leq N_d = 
GW_0 (X , d , pt, \dots , pt),$$
for every generic $J \in \R {\cal J}_\omega$ and every $\underline{x} \in X^{k_d}$ such that
$\underline{x} \cap \R X = r$.
\end{cor}
\begin{rem}
\label{remcomment}
Let us end this paragraph with some comments on the hypothesis of Theorem \ref{theosemipos}.

1) The condition $2 / (c_1(X)d - 2)$ is a necessary condition for the genus zero Gromov-Witten 
invariant $GW_0 (X , d , pt, \dots , pt)$ to be non-trivial. This just comes from dimensional reasons
since we do consider only point conditions throughout this paper.

2) The strongly semipositive condition as well as the condition $c_1 (X)d > 2$ is made to prevent
the appearance of non simple real stable maps in the Gromov compactification of the moduli space
of real rational pseudo-holomorphic curves over a generic path in $\R {\cal J}_\omega$. The 
treatment of such
non simple stable maps requires more involved techniques, see \cite{LT} for example. Note that the
theory of polyfolds under construction might provide helpful techniques to remove these assumptions,
see \cite{polyfold}.

3) The hypothesis HYP are topological contitions on degree two characteristic classes of the real locus
$\R X$. They are used to define the spinor state of real rational curves which are crucial in the
defition of the invariant $\chi^{\p}$. Can similar invariants be obtained without these topological 
contitions? I don't know. Note that similar issues appeared in \cite{FOOO} in order to
prove the orientability of the moduli space of pseudo-holomorphic discs having boundary in some
Lagrangian submanifold $L$. There, $L$ was assumed to be relatively spin.

4) The existence of at least one real point in the configuration is to prevent
the appearance of real rational curves with empty real locus in ${\cal R}_d (\underline{x} , J , h)$,
see Lemma \ref{lemmarealpart}. For such real rational curves with empty real locus, spinor states are
not defined and I cannot obtain similar invariants yet. This subtle problem appears in
the important case of Calabi-Yau threefolds, where $\underline{x} = \emptyset$. I believe however
that there should be a way to overcome sometimes this difficulty, but cannot do it yet.
\end{rem}

\subsection{Topological interpretation}

Note that the singularities of $\R_\tau \overline{\cal M}^d_{k_d} (X)$ are of codimension two at 
least since
it is normal, see Theorem $2$ of \cite{FP}. This space thus carries a first Stiefel-Whitney 
class. For $D \in H_{3k_d - 1} (\R_\tau \overline{\cal M}^d_{k_d} (X) ; \Z/2\Z)$, denote by
$D^\vee$ its image under the morphism $H_{3k_d - 1} (\R_\tau \overline{\cal M}^d_{k_d} (X) ; \Z/2\Z)
\to H^1 (\R_\tau \overline{\cal M}^d_{k_d} (X) ; \Z/2\Z)$.
\begin{prop}
\label{propw1}
The first Stiefel-Whitney class of every component $\R {\cal M}^*$ of $\R_\tau {\cal M}^d_{k_d} (X)$ 
which contains a balanced curve writes
$$w_1 (\R {\cal M}^*) = (\R_\tau ev^d)^* w_1 (\R_\tau X^{k_d}) + 
\sum_{D \subset \text{Red}'} \epsilon (D) D^\vee \in H^1 (\R {\cal M}^* ; \Z/2\Z),$$
where $\epsilon (D) \in \{ 0, 1 \}$ and if $\epsilon (D) = 1$, the irreducible component $D$ of 
$\text{Red}$ gets contracted by the evaluation map $\R_\tau ev^d$. $\square$
\end{prop}
Here $\text{Red}'$ denotes the union of the divisor $\text{Red}$ introduced in Theorem \ref{theoprelim}
and the divisors of non-balanced curves $(u , C , \underline{z})$ such that $\dim H^1 (C ; N_u \otimes
{\cal O}_C (-  \underline{z})) \geq 2$, in case such divisors exist. Denote by $\text{Red}_1$ the 
union of the irreducible components $D$ of $\text{Red}'$ for which $\epsilon (D) = 1$.
Equip $\R_\tau X^{k_d}$ with the system of twisted integer coefficients ${\cal Z}$ and denote
by $[\R_\tau X^{k_d}] \in H_{3k_d} (\R_\tau X^{k_d} ; {\cal Z})$ one associated fundamental class. 
Denote by ${\cal Z}^*$ the local system of coefficients of 
$\R {\cal M}^*$ pulled back from ${\cal Z}$ by $\R_\tau ev^d$,
see \cite{Steen}. 
\begin{prop}
\label{propfond}
Under the hypothesis of Proposition \ref{propw1}, there exists a unique fundamental class 
$[\R {\cal M}^*] \in H_{3k_d} (\R {\cal M}^* , \text{Red}_1 ; {\cal Z}^*)$ such that for 
every balanced curve $(u, C , \underline{z}) \in \R {\cal M}^*$, the morphism
$(\R_\tau ev^d)_* : H_{3k_d} (\R {\cal M}^* , 
\R {\cal M}^* \setminus \{ (u, C , \underline{z}) \} ; {\cal Z}^*) \to
H_{3k_d} (\R_\tau X^{k_d} , \R_\tau X^{k_d} \setminus \{ u(\underline{z}) \} ; {\cal Z})$ sends
$[\R {\cal M}^*]$ onto $sp_{\p , {\goth o}} (u, C , \underline{z}) [\R_\tau X^{k_d}]$. $\square$
\end{prop}
Since $\R_\tau ev^d (\text{Red}_1)$ is of codimension two, the group $H_{3k_d} (\R_\tau X^{k_d} , 
\R_\tau ev^d (\text{Red}_1) ; {\cal Z})$ is cyclic, generated by $[\R_\tau X^{k_d}]$. The integer 
$\chi_r^{d, h \p}$ is nothing but the one defined by the relation $(\R_\tau ev^d)_* 
[\R_\tau \overline{\cal M}^d_{k_d} (X)] = \chi_r^{d, h \p} [\R_\tau X^{k_d}]$, where the fundamental 
class $[\R_\tau \overline{\cal M}^d_{k_d} (X)]$ is given by Proposition \ref{propfond} and restricted 
to connected components $\R {\cal M}^*$ for which $d_{\R {\cal M}^*} = h$, see Lemma \ref{lemmarm*}.

\section{Proof of Theorem \ref{theosemipos}}

Let $J_0 , J_1 \in \R {\cal J}_\omega$ be two generic real almost complex structures, so that 
the integers $\chi_r^{d, h \p , {\goth o}} (\underline{x} , J_0)$ and 
$\chi_r^{d, h \p , {\goth o}} (\underline{x} , J_1)$ are well defined. We have to prove that these 
integers are the same. Let $\gamma : t \in [0,1] \mapsto J_t \in \R {\cal J}_\omega$ be a generic 
path joining $J_0$ to $J_1$, transversal to $\pi_\R$. Let $\R {\cal M}_\gamma = 
\R {\cal M}^d (\underline{x})^* \times_\gamma [0,1]$, $\R \overline{\cal M}_\gamma$ be its Gromov 
compactification and $\pi_\gamma : \R \overline{\cal M}_\gamma \to [0,1]$ be the associated projection.
Then $\R \overline{\cal M}_\gamma$ provides a cobordism between ${\cal R}_d (\underline{x} , J_0 , h)$
and ${\cal R}_d (\underline{x} , J_1 , h)$. On each connected component of the complement
of both the the critical values of $\pi_\gamma$ and the elements of
$\pi_\gamma (\R \overline{\cal M}_\gamma \setminus \R {\cal M}_\gamma)$, the integer 
$\chi_r^{d, h \p , {\goth o}} (\underline{x} , J_t)$ is constant. We thus just have to prove that
this integer also does not change while crossing one of these values. Theorem 
\ref{theosemipos} hence follows from Theorems \ref{theocritsp} and \ref{theoredsp}. $\square$

\subsection{Generic critical points of $\pi_\R$}

\begin{lemma}
\label{lemmacusp}
The space of elements $[u, J_S , J] \in \R {\cal M}^d (\underline{x})^*$ for which $u$ is not an 
immersion is a sub stratum of codimension $n-1$ of $\R {\cal M}^d (\underline{x})^*$.
\end{lemma}

{\bf Proof:}

This fact was proved in the first part of Proposition $2.7$ of \cite{Wels1} when $n=2$ -and mainly
follows from the results of \S $3$ of \cite{Shev}-. The proof being readily the same in higher 
dimensions, it is not reproduced here. $\square$ \\

Let $\gamma : t \in [0,1] \mapsto J_t \in \R {\cal J}_\omega$ be a generic path transversal to
$\pi_\R$. Let $\R {\cal M}_\gamma = \R {\cal M}^d (\underline{x})^* \times_\gamma [0,1]$ and
$\pi_\gamma : \R {\cal M}_\gamma \to [0,1]$ be the associated projection. From Lemma \ref{lemmacusp}
follows that all the elements of $\R {\cal M}_\gamma$ are immersions.
\begin{theo}
\label{theocritsp}
Let $[u_{t_0}, J_S^{t_0} , J_{t_0}] \in \R {\cal M}_\gamma$ be a critical point of $\pi_\gamma$. 
Let $\mu : \lambda \in ] - \epsilon , \epsilon [ \mapsto \mu (\lambda) \in \R {\cal M}_\gamma$ be a 
local parameterization such that $\mu (0) =  [u_{t_0}, J_S^{t_0} , J_{t_0}]$. Then, as soon as
$\epsilon$ is small enough, the following alternative holds. Either $\pi_\gamma \circ \mu$ crosses
$t_0$ as $\lambda$ crosses $0$ and then $sp_{\p , {\goth o}} (\mu (\lambda))$ does not depend on
$\lambda \in ] - \epsilon , \epsilon [ \setminus \{ 0 \}$, or $\pi_\gamma \circ \mu$ does not cross
$t_0$ as $\lambda$ crosses $0$ and then $sp_{\p , {\goth o}} (\mu (\lambda)) = -
sp_{\p , {\goth o}} (\mu (- \lambda))$ for every 
$\lambda \in ] - \epsilon , \epsilon [ \setminus \{ 0 \}$.
\end{theo}

{\bf Proof:}

Denote by $\mu (\lambda) = [u_\lambda, J_S^\lambda , J_\lambda]$ and fix a $\Z / 2\Z$-equivariant
trivialization $N \to S$ of the complex normal bundles $N_{u_\lambda} \to S$. We deduce a
family of generalized real Cauchy-Riemann operators $D_\lambda : L^{k,p} (S , N_{-\underline{z}})_{+1}
\to L^{k-1,p} (S , \Lambda^{0,1} S \otimes N)_{+1}$ parameterized by $\lambda \in ] - \epsilon , 
\epsilon [$. Let $E_0$ be a closed complement to the one dimensional kernel of $D_0$ and
for $\lambda \in ] - \epsilon , \epsilon [$, $F_\lambda = D_\lambda (E_0)$. Then, 
$D_\lambda$ induces a family of morphisms $\widetilde{D}_\lambda : H^0 = 
L^{k,p} (S , N_{-\underline{z}})_{+1} / E_0 \to H^1_\lambda = 
L^{k-1,p} (S , \Lambda^{0,1} S \otimes N)_{+1} / F_\lambda$. Note that the one dimensional vector
space $H^0$ (resp. $H^1_\lambda$) is trivialized by $v_\lambda = \frac{d}{d \lambda} u_\lambda$
(resp. $\frac{d}{d t} \gamma|_{t = \pi_\gamma \circ \mu (\lambda)}$), so that the linear maps
$\widetilde{D}_\lambda$ and $\frac{d}{d \lambda} (\pi_\gamma \circ \mu )$ get conjugated.
More generaly, for every operator $D \in \R {\cal O}p_{\overline{\partial} + R} (N)$, denote by
$\widetilde{D}$ the induced morphism $L^{k,p} (S , N_{-\underline{z}})_{+1} / E_0 \to
L^{k-1,p} (S , \Lambda^{0,1} S \otimes N)_{+1} / D(E_0)$. The hypersurface ${\cal H}^1$ of
$\R {\cal O}p_{\overline{\partial} + R} (N)$ made of operators having a one dimensional cokernel
is defined in a neighbourhood of $D_0$ as ${\cal H}^1 = \{ D \in \R {\cal O}p_{\overline{\partial} 
+ R} (N) \, | \, \widetilde{D} = 0 \}$, see Lemma \ref{lemmawall}. Thus, the curve
$\lambda \in ] - \epsilon , \epsilon [ \mapsto D_\lambda \in 
\R {\cal O}p_{\overline{\partial} + R} (N)$ crosses ${\cal H}^1$ at $\lambda = 0$ if and only if
$\frac{d}{d \lambda} (\pi_\gamma \circ \mu ) (v_\lambda)$ crosses $t_0$ as $\lambda$ crosses $0$
that is if and only if $\pi_\gamma$ has a local extremum at $\lambda = 0$. Now from Proposition 
\ref{propwall}, if $D_\lambda$ crosses ${\cal H}^1$ at $\lambda = 0$,
then $sp_{\p , {\goth o}} (D_\lambda) = - sp_{\p , {\goth o}} (D_{- \lambda})$,
$\lambda \in ] - \epsilon , \epsilon [ \setminus \{ 0 \}$, whereas $sp_{\p , {\goth o}} (D_\lambda) 
= sp_{\p , {\goth o}} (D_{- \lambda})$ otherwise. $\square$

\subsection{Passing through real reducible curves}

Let $\gamma : t \in [0,1] \mapsto J_t \in \R {\cal J}_\omega$ be a generic path transversal to
$\pi_\R$. Let $\R {\cal M}_\gamma = \R {\cal M}^d (\underline{x})^* \times_\gamma [0,1]$,
$\R \overline{\cal M}_\gamma$ be its Gromov compactification and
$\pi_\gamma : \R \overline{\cal M}_\gamma \to [0,1]$ be the associated projection.
\begin{lemma}
\label{lemmaaltern}
Under the assumptions of Theorem \ref{theosemipos}, as soon as $\gamma$ is generic enough,
$\R \overline{\cal M}_\gamma \setminus \R {\cal M}_\gamma$ consists of finitely many reducible curves
having two irreducible components, both real and embedded, intersecting in a unique real ordinary 
double point. Moreover, if $d_1 , d_2$ (resp. $k_1 , k_2$) denote the homology classes of these
components (resp. the number of marked points on them), so that $d_1 + d_2 = d$ and $k_1 + k_2 = k_d$,
then either $c_1 (X)d_1 - 1 = 2k_1$ and $c_1 (X)d_2 - 1 = 2(k_2 - 1)$, or
$k_i = E \big( \frac{1}{2} (c_1 (X)d_i - 2) \big) + 1$, where $E()$ denotes the integer part.
\end{lemma}

{\bf Proof:}

From Gromov compactness theorem, elements of $\R \overline{\cal M}_\gamma \setminus 
\R {\cal M}_\gamma$ correspond to reducible curves parameterized by a tree of complex spheres.
Since the manifold is semipositive, the moduli space of such curves comes with a Fredholm projection
onto ${\cal J}_\omega$ whose Fredholm index is one minus the number of spheres in this tree,
see Theorem $6.2.6$ of \cite{McDSal}. Note that since
marked points are not at singular points of the parameterizing tree and the configuration
of points $\underline{x}$ contains at least one real point, both irreducible components of the
reducible curve should be real. The first part of the lemma follows. Now the numerical conditions
on $k_1 , k_2$ are obtained exactly in the same way as in Proposition $3.3$ of \cite{Wels2}, it is not reproduced here. $\square$\\

Denote by $\R {\cal M}^{d_1 , d_2}_{k_1 , k_2} (\underline{x})$ the universal moduli space of
simple real stable maps having two irreducible components $C_1$, $C_2$, both real with $k_1 , k_2$
marked points on them respectively and which realize the homology classes $d_1 , d_2$. Denote by
$\pi_\R^{d_1 , d_2}$ the index $-1$ Fredholm projection 
$\R {\cal M}^{d_1 , d_2}_{k_1 , k_2} (\underline{x}) \to \R {\cal J}_\omega$. 
\begin{prop}
\label{propint}
Let $(u, J_S, J) \in \R {\cal M}^{d_1 , d_2}_{k_1 , k_2} (\underline{x})$ be given by Lemma
\ref{lemmaaltern}. Then, there exists a path $(J_\lambda)_{\lambda \in [0,1]}$ in
$\R {\cal J}_\omega$ such that $J_0 = J$, $(u, J_S, J_\lambda) \in 
\R {\cal M}^{d_1 , d_2}_{k_1 , k_2} (\underline{x})$ for every $\lambda \in [0,1]$ and
$J_1$ is integrable in a neighbourhood of $u(C_1 \cup C_2)$. Moreover, $J_1$ can be chosen in
the form given by Propositions $1.4$ or $1.6$ of  \cite{Wels2} depending on whether
$k_1 = E \big( \frac{1}{2} (c_1 (X)d - 2) \big)$ or 
$k_1 = E \big( \frac{1}{2} (c_1 (X)d - 2) \big) + 1$.
\end{prop}

{\bf Proof:}

Such an homotopy $(J_\lambda)_{\lambda \in [0,1]}$ can be obtained in the following way. We first
stretch the almost complex structure $J$ in a neighbourhood of the double point using a one
parameter family of $\Z / 2\Z$-equivariant diffeomorphisms of $X$ which read as a family of homotheties
in a local chart mapping the curve onto two coordinate axis of $\C^3$ and mapping $J$ at the singular
point onto the complex structure of $\C^3$. Having the scale of the homothety converging to $+\infty$,
we deduce a homotopy $(J_\lambda)_{\lambda \in [0,\frac{1}{2}]}$ such that $J_\frac{1}{2}$ is
integrable in a neighbourhood of the singular point. For every $i \in \{ 1 , 2 \}$, fix a
$\Z / 2\Z$-equivariant identification of a tubular neighbourhood of $u (C_i)$ in $X$ with a
neighbourhood of the zero section in its normal bundle $N_i = u^*TX / u_* TC_i$. The latter is equipped
with the almost complex structure $J_\frac{1}{2}$. We then stretch $J_\frac{1}{2}$ using
a one parameter family of homotheties in the fibres of $N_i$ whose scale converge to $+\infty$. The
path of almost complex structures we obtain does converge since the zero section is pseudo-holomorphic
and we get a homotopy $(J_\lambda)_{\lambda \in [\frac{1}{2} , \frac{3}{4}]}$ such that
$J_\frac{3}{4}$ is integrable in a neighbourhood of the singular point and equip $N_i$ with the
structure of a complex vector bundle. Now, let $J_1$ be a holomorphic structure on this complex vector
bundle and $\nabla_1$ be an associated complex connection. Then, for every $y \in N_i$, 
$J_\frac{3}{4} (y)$ writes $J_1 (y) + A(y)$ where $A(y) \in \Lambda^{0,1} C_i \otimes N_i$.
We thus set for $\lambda \in [ \frac{3}{4} , 1]$ and $y \in N_i$, $J_{\lambda} (y) =
J_1 (y) + 4(1 - \lambda) A(y)$. The result follows from the fact that $J_1$ can be chosen in
the form given by Propositions $1.4$ or $1.6$ of  \cite{Wels2} depending on whether
$k_1 = E \big( \frac{1}{2} (c_1 (X)d - 2) \big)$ or 
$k_1 = E \big( \frac{1}{2} (c_1 (X)d - 2) \big) + 1$. $\square$
\begin{theo}
\label{theoredsp}
Let $\gamma : t \in [0,1] \mapsto J_t \in \R {\cal J}_\omega$ be a generic path transversal to
$\pi_\R$. Let $\R {\cal M}_\gamma = \R {\cal M}^d (\underline{x})^* \times_\gamma [0,1]$,
$\R \overline{\cal M}_\gamma$ be its Gromov compactification and
$\pi_\gamma : \R \overline{\cal M}_\gamma \to [0,1]$ be the associated projection.
Let $(u, J_S, J) \in \R \overline{\cal M}_\gamma \setminus \R {\cal M}_\gamma$ be given by Lemma
\ref{lemmaaltern} and $t_0 = \pi_\gamma (u, J_S, J) \in ] 0 , 1 [$. Then, there exist a neighbourhood
$W$ of $(u, J_S, J)$ in $\R \overline{\cal M}_\gamma$ and $\epsilon > 0$ such that
$\sum_{C \in \pi_\gamma^{-1} (t) \cap W} sp_{\p , {\goth o}} (C)$ does not depend on $t \in
]t_0 - \epsilon , t_0 + \epsilon [ \setminus \{ t_0 \}$.
\end{theo}

{\bf Proof:}

From Lemma \ref{lemmaaltern}, $(u, J_S, J) \in \R {\cal M}^{d_1 , d_2}_{k_1 , k_2} (\underline{x})$
with either $k_1 = E \big( \frac{1}{2} (c_1 (X)d - 2) \big)$ or 
$k_1 = E \big( \frac{1}{2} (c_1 (X)d - 2) \big) + 1$. From Proposition \ref{propint}, we can assume 
that $J$ is integrable in a neighbourhood of $u(C_1 \cup C_2)$ and of the form given by Propositions 
$1.4$ or $1.6$ of  \cite{Wels2} depending on whether
$k_1 = E \big( \frac{1}{2} (c_1 (X)d - 2) \big)$ or 
$k_1 = E \big( \frac{1}{2} (c_1 (X)d - 2) \big) + 1$. Indeed, let $\mu : \lambda \in [0,1] \to
(u, J_S, J_\lambda) \in \R {\cal M}^{d_1 , d_2}_{k_1 , k_2} (\underline{x})$ be the path given by
this Proposition  \ref{propint}. Then, after perturbing $\mu$ if necessary, we can assume that
it crosses transversely the critical locus of $\pi_\R^{d_1 , d_2} : 
\R {\cal M}^{d_1 , d_2}_{k_1 , k_2} (\underline{x}) \to \R {\cal J}_\omega$ at finitely many points,
where the cokernel of  $\pi_\R^{d_1 , d_2}$ is two dimensional. Outside of these points,
$\pi_\R^{d_1 , d_2}|_{\text{Im} (\mu)}$ is an immersion and thus $\text{Im} (\pi_\R^{d_1 , d_2})$
is locally a wall which divides $\R {\cal J}_\omega$ in two connected components. Let ${\cal W}$
be a neighbourhood of $\text{Im} (\mu)$ in $\R \overline{\cal M}^d (\underline{x})$. Then, from
Theorem \ref{theocritsp}, as soon as ${\cal W}$ is small enough, 
$\sum_{C \in \pi_\R^{-1} (J) \cap {\cal W}} sp_{\p , {\goth o}} (C)$ takes one value on each side of
the wall $\text{Im} (\pi_\R^{d_1 , d_2})$. Now let $(u, J_S, J_\frac{1}{2}) \in 
\R {\cal M}^{d_1 , d_2}_{k_1 , k_2} (\underline{x})$ be a critical point of $\pi_\R^{d_1 , d_2}$
and $J_\frac{1}{2} \in D \subset \R {\cal J}_\omega$  be a closed disc transversal to
$\pi_\R^{d_1 , d_2}$ at $J_\frac{1}{2}$. Then, $(\pi_\R^{d_1 , d_2})^{-1} (D)$ is a smooth curve
of $\R {\cal M}^{d_1 , d_2}_{k_1 , k_2} (\underline{x})$ which projects onto a curve of $D$
which is cuspidal at $J_\frac{1}{2}$. The connected component of this curve which contains
$(u, J_S, J_\frac{1}{2})$ is homeomorphic to an interval whose image curve $B$ intersects $\partial D$
at two points. The complement of these two points in $\partial D$ consists of two intervals
$\partial D_1$ and $\partial D_2$ and once more, from
Theorem \ref{theocritsp}, as soon as ${\cal W}$ is small enough, the value 
$\sum_{C \in \pi_\R^{-1} (J) \cap {\cal W}} sp_{\p , {\goth o}} (C)$ is the same on each of these
intervals. 
$${\vcenter
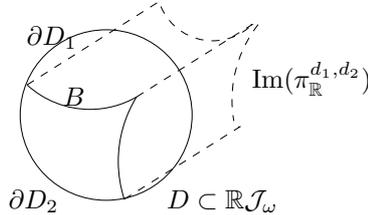
{\hbox{\begin{picture}(0,0)%
\includegraphics{semi1.pstex}%
\end{picture}%
\setlength{\unitlength}{3315sp}%
\begingroup\makeatletter\ifx\SetFigFont\undefined%
\gdef\SetFigFont#1#2#3#4#5{%
  \reset@font\fontsize{#1}{#2pt}%
  \fontfamily{#3}\fontseries{#4}\fontshape{#5}%
  \selectfont}%
\fi\endgroup%
\begin{picture}(1853,1600)(2206,-4619)
\put(2611,-3796){\makebox(0,0)[lb]{\smash{\SetFigFont{10}{12.0}{\rmdefault}{\mddefault}{\updefault}{$B$}%
}}}
\put(2206,-4561){\makebox(0,0)[lb]{\smash{\SetFigFont{10}{12.0}{\rmdefault}{\mddefault}{\updefault}{$\partial D_2$}%
}}}
\put(2341,-3346){\makebox(0,0)[lb]{\smash{\SetFigFont{10}{12.0}{\rmdefault}{\mddefault}{\updefault}{$\partial D_1$}%
}}}
\put(3376,-4561){\makebox(0,0)[lb]{\smash{\SetFigFont{10}{12.0}{\rmdefault}{\mddefault}{\updefault}{$D \subset \R {\cal J}_\omega$}%
}}}
\put(4006,-3661){\makebox(0,0)[lb]{\smash{\SetFigFont{10}{12.0}{\rmdefault}{\mddefault}{\updefault}{$\text{Im} (\pi_\R^{d_1 , d_2})$}%
}}}
\end{picture}
}}}$$
Thus, the value $\sum_{C \in \pi_\R^{-1} (J) \cap {\cal W}} sp_{\p , {\goth o}} (C)$ on
each side of the wall $\text{Im} (\pi_\R^{d_1 , d_2})$ is the same at $J_0$ and $J_1$. To prove
Theorem \ref{theoredsp}, we then just have to prove that these values are the same on both sides of
the wall at $J_1$. Now, to cross this wall, we can fix $J$ and have one real point of the
configuration $\underline{x}$ vary since this can be performed equivalently by fixing $\underline{x}$
and producing some Hamiltonian deformation of $J$. The end of the proof thus goes now along the same
lines as the one of \cite{Wels2} since all the arguments used there where local. $\square$

\begin{rem}
The difference between this new version and the original version of our preprint is that the main result Theorem \ref{theosemipos}
is restricted
to dimension six. This is because our Theorem $3.4$ of the original version does not seem correct in general in dimension 
greater than six, as pointed out by a referee. Note that in the meanwhile, these invariants have been interpreted using the work 
\cite{FOOO}
by Cheol-Hyun Cho in \cite{Cho} and Jake Solomon in \cite{Sol}. In \cite{Sol}, these invariants are extended to Calabi-Yau 
six-manifolds through a treatment of multiple curves than can appear (see Remark \ref{remcomment}) and computed in
\cite{PSW} for real quintic three-folds. Also, the exact value
of $\epsilon (D)$ in Proposition \ref{propw1} has been computed by Nicolas Puignau in \cite{Pui} in the case of the projective plane
or quadric hyperboloid. Finally, sharpness of the lower bounds in Corollary
\ref{corlowerbounds} as well as some computations have been obtained when $r=1$ for 
the quadric ellipsoid three-fold in \cite{WelsSFT}.
\end{rem}

\addcontentsline{toc}{part}{\hspace*{\indentation}Bibliography}


\bibliography{semipositive}
\bibliographystyle{abbrv}

\noindent \'Ecole normale sup\'erieure de Lyon\\
Unit\'e de math\'ematiques pures et appliqu\'ees\\
UMR CNRS $5669$\\
$46$, all\'ee d'Italie\\
$69364$, Lyon cedex $07$\\
(FRANCE)\\
e-mail : {\tt jwelschi@umpa.ens-lyon.fr}

\end{document}